\documentclass{article}
\usepackage{amssymb}
\usepackage{epsfig}
\usepackage{amsmath}
\usepackage{amsfonts}
\usepackage{tikz}
\usepackage{float}
\usepackage{color}
\usepackage{hyperref}

\newcommand{\R}{{{\Bbb R}}}

\newtheorem{theorem}{\sc Theorem}[section]
\newtheorem{proposition}{\sc Proposition}[section]
\newtheorem{lemma}{\sc Lemma}[section]
\newtheorem{definition}{\sc Definition}[section]
\newtheorem{remark}{\sc Remark}[section]
\newtheorem{notation}{\sc Notation}[section]
\newtheorem{corollary}{\sc Corollary}[section]
\newtheorem{example}{\sc Example}[section]

\usepackage{fancyheadings}
\renewcommand{\sectionmark}[1]%
                    {\markboth{#1}{}}
\renewcommand{\sectionmark}[1]%
                    {\markright{\thesection\ #1}}
\title{
On some non-linear boundary value problems related to a Black--Scholes model with transaction costs}

\author{Rub\'en Figueroa and Maria do Ros\'ario Grossinho}

\date{}

\begin{document}
\maketitle


\def\tanh#1{\,{\normalsize tanh}{\,#1}\,}
\def\simbolo#1#2{#1\dotfill{#2}}
\def\qed{\hbox to 0pt{}\hfill$\rlap{$\sqcap$}\sqcup$\medbreak}
\def\theequation{\arabic{section}.\arabic{equation}}
\def\thesection {\arabic{section}}

\begin{abstract} We deal with some generalizations on a Black--Scholes model arising in financial mathematics. As novelty in this paper, we consider a variable volatility and abstract functional boundary conditions, which allow us to treat a very large class of problems involving Black--Scholes equation. Our main results involve the existence of extremal solutions in presence of lower and upper solutions. Some examples of application are provided too.
\end{abstract}

\section{Introduction}

In this paper we are concerned with the following non-linear boundary value problem

\begin{equation}\label{p}
\left\{
\begin{array}{ll}
x^3 (V'')^2(x) + p(x) x^2 V''(x) + q(x) (xV'(x)-V(x))=0, \\
\\
B_1(V(c),V)=0, \quad B_2(V(d),V)=0,
\end{array}
\right.
\end{equation}
where $p,q$ are nonnegative bounded functions which could be discontinuous in $[c,d]$, $c,d >0$, and $B_i:\R \times \mathcal{C}([c,d]) \longrightarrow \R$, $i=1,2$, are functions which satisfy some conditions that we will state later. We observe that under this framework, a large class of boundary conditions is included, namely:

\begin{enumerate}

\item{Dirichlet conditions: $B_1(V(c),V)=V(c)-V_c, \quad B_2(V(d),V)=V(d)-V_d$,}

\item{Initial--integral conditions: $B_1(V(c),V)=V(c)-\displaystyle{\int_c^d k(x) V(x) \, dx}$,}

\item{Multipoint conditions: $B_2(V(d),V)=V(d)-\displaystyle{\sum_{j=1}^n V(x_j)}$,}

\end{enumerate}

This study follows and generalizes the results contained in \cite{gromor} with respect to the problem 
\begin{equation}\label{bs3}
\left\{
\begin{array}{ll}
x^3 (V'')^2(x) + p\, x^2 V''(x) + q (x \, V'(x) - V(x))=0, \\
\\
V(c)=V_c,\quad V(d)=V_d.
\end{array}
\right.
\end{equation}
In \cite{gromor}, it is assumed that $p,q$ are positive constants and  $V_c < V_d$. The contributions of the present paper are the following.
First we address this problem but drop the condition $V_c < V_d$, and replace the constants $p,q$ by two functions $p(x),q(x)$; second, we replace the Dirichlet conditions by functional boundary conditions, which allows us to consider a very large class of problems for the equation of (\ref{bs3}).

These problems are related to financial option pricing, since they address the existence of stationary solutions of a class of generalizations of the classical Black-Scholes model (BS), introduced in 1973 \cite{blsc}, with equation:
\begin{equation}\label{bs1}
\dfrac{\partial V}{\partial t} + \dfrac{1}{2} \sigma^2 S^2 \dfrac{\partial^2 V}{\partial S^2} + r\left( S \dfrac{\partial V}{\partial S}-V \right)=0,
\end{equation}
where $V$ represents the value of a call or put option, depending  on an underlying asset $S$ and on time $t$, $r$ is the interest short rate and $\sigma$ is the volatility of the asset price. In the (BS) model, $S$  is modelled as a geometric Brownian motion and  no costs are considered when financial transactions hold. 

Suppose that transaction costs are included in the model under the assumption that they are a percentage of the transaction, given as in \cite{amsab} by a linear
function $h$\ of the number of shares traded, i.e.,
   $h(\eta )=a-b \lambda_2,$
  where $\eta $\ is the number of shares traded and $a, b>0.$ Then, if $a$ is small enough and $V_{SS}>0$ (see \cite{amsab}, \cite{gromor2}, \cite{gromor}) the following nonlinear version of (\ref{bs1}) is obtained:

\begin{equation}\label{bs2}
\dfrac{\partial V}{\partial t} + \dfrac{1}{2} \tilde{\sigma}^2 S^2 \dfrac{\partial^2 V}{\partial S^2} + b\sigma^2 S^3 \left(\dfrac{\partial^2 V}{\partial S^2}\right)^2+r\left(S \dfrac{\partial V}{\partial S}-V \right)=0,
\end{equation}
where $\tilde{\sigma}$ is an adjusted volatility. Now, if we consider the stationary version of (\ref{bs2}), we obtain the above ordinary differential equation (\ref{bs3})
where $p=\dfrac{\tilde{\sigma}^2}{2b \sigma^2}$ and $q=\dfrac{r}{b \sigma^2}$  are constants.

This paper is organized as follows: in Section $2$, we introduce an auxiliary non-linear boundary value problem and the notions of upper and lower solutions used later. In Section $3$, we consider the problem (\ref{p}) with Di\-ri\-chlet boundary conditions. So, we start from paper \cite{gromor} and we generalize it by considering $p(x),q(x)$ instead of constants $p,q$, which corresponds to variable volatility in the (BS) model, and by dropping the condition $V_c < V_d$. In Section $4$, we deal with the problem (\ref{p}) on its full version, that is, with functional boundary conditions. Namely, we provide a result on the existence of extremal solutions between lower and upper solutions, by using a generalized iteration with Dirichlet problems. Some examples of application are provided, too.

\section{Auxiliary problem and upper and lower solutions}

If we look to the equation of (\ref{p})
\begin{equation*}
x^3 (V'')^2(x) + p(x) S^2 V''(x) + q(x) (x \, V'(x) - V(x))=0,
\end{equation*}
and, as in \cite{amsab}, \cite{gromor2} and \cite{gromor}, solve it algebraically  in order to $V''$, we obtain the following equation 
\begin{equation}\label{v}
V''(x)+H(x,V(x),V'(x))=0,
\end{equation}
where
\begin{equation}\label{HH}
H(x,y,z)=\dfrac{p(x)x^2 - \sqrt{p^2(x) x^4 + 4x^3 q(x)|xz-y|}}{2x^3}.
\end{equation}
So, we will study the problem (\ref{p}) by considering first the following related problem:

\begin{equation}\label{pp}
\left\{
\begin{array}{ll}
V''(x)+H(x,V(x),V'(x))=0, \\
\\
B_1(V(c),V)=0, \quad B_2(V(d),V)=0,
\end{array}
\right.
\end{equation}
and then we will see how solutions of this problem are solutions of our original problem.

We will use the method of upper and lower solutions for this problem and we will begin by considering the classical notions of $\mathcal{C}^2-$lower and upper solutions. However, in Sections $3$ and $4$, we will use some weaker notions since we will need lower and upper solutions to have ``angles''. So, denoting by $D_{-}f(x)$ and $D^{+}f(x)$, respectively, the lower left-hand and the upper right-hand Dini-derivatives of the function $f$ at $x$, we introduce the following definitions (see \cite{coha}).

\begin{definition}\label{lu2} We say that $\alpha \in AC([c,d])$ is a lower solution for problem (\ref{pp}) if
$$
B_1(\alpha(c),\alpha) \le 0, \quad B_2(\alpha(d),\alpha) \le 0,$$
and for each $x_0 \in (c,d)$ one of the following conditions holds:
\begin{enumerate}

\item{$D_-\alpha(x_0) < D^+\alpha(x_0)$;}

\item{There exists an open interval $J_0$ such that $x_0 \in J_0 \subset (c,d)$, $\alpha_{|J_0} \in W^{2,1}(J_0)$ and for almost all $x \in J_0$ we have
$$
\alpha''(x)+H(x,\alpha(x),\alpha'(x)) \ge 0.$$}

\end{enumerate}

We say that $\beta \in AC([c,d])$ is an upper solution for problem (\ref{pp}) if
$$
B_1(\beta(c),\beta) \ge 0, \quad B_2(\beta(d),\beta) \ge 0,$$
and for each $x_0 \in (c,d)$ one of the following conditions holds:
\begin{enumerate}

\item{$D^-\beta(x_0) > D_+\beta(x_0)$;}

\item{There exists an open interval $J_0$ such that $x_0 \in J_0 \subset (c,d)$, $\beta_{|J_0} \in W^{2,1}(J_0)$ and for almost all $x \in J_0$ we have
$$
\beta''(x)+H(x,\beta(x),\beta'(x)) \le 0.$$}

\end{enumerate}
\end{definition}
Notice that if $\alpha$, $\beta$ are classical $\mathcal{C}^2-$ lower and upper solutions for problem (\ref{pp}) then they are also lower and upper solutions in the sense referred above.

\section{Non-linear problem with Dirichlet conditions}

Consider the problem (\ref{p}) with $$
B_1(V(c),V)=V(c)-V_c, \quad B_2(V(d),V)=V(d)-V_d$$ that is, with 
standard Dirichlet conditions. So, we have in this case
\begin{equation}\label{p2}
\left\{
\begin{array}{ll}
x^3 (V'')^2(x) + p(x) x^2 V''(x) + q(x) (xV'(x)-V(x))=0, \\
\\
V(c)=V_c, \quad V(d)=V_d.
\end{array}
\right.
\end{equation}
The auxiliary problem referred in the previous sections is now
\begin{equation}\label{p1}
\left\{
\begin{array}{ll}
V''(x)+H(x,V(x),V'(x))=0, \\
\\
V(c)=V_c, \quad V(d)=V_d,
\end{array}
\right .
\end{equation}
where the function $H$ is given by (\ref{HH}).

From the study of the problem (\ref{p1}), we will deduce later existence and localization results for the problem (\ref{p2}). 
Next proposition establishes adequately the existence of classical $\mathcal{C}^2-$ upper and lower solutions for (\ref{p1})

\begin{proposition} The following assertions hold:

\begin{enumerate}

\item[$(1)$]{If  $\dfrac{V_d}{d} \le \dfrac{V_c}{c}$, then the function
\begin{equation}\label{alpha1}
\alpha_1(x)=\dfrac{V_d}{d} x
\end{equation}
is a $\mathcal{C}^2$-lower solution for the problem (\ref{p1}).

The converse  is also valid.}
\item[$(2)$]{Take $k >0$ such that
\begin{equation}\label{k}
k \ge \sqrt{\dfrac{Q}{c^3}} \sqrt{ \max_{x \in [c,d]} \left|\dfrac{k}{2}(x^2-cd)+  \dfrac{cV_d -dV_c }{d-c}\right|},
\end{equation}
where
$$
Q=\max_{x \in [c,d]} q(x).$$
Then the function
\begin{equation}\label{alpha2}
\alpha_k(x)=\dfrac{k}{2} x^2 + \left(\dfrac{V_d-V_c}{d-c}-\dfrac{k}{2}(d+c)\right)x + \dfrac{k}{2} cd - \dfrac{cV_d -dV_c }{d-c}
\end{equation}
is a $\mathcal{C}^2$-lower solution for the problem (\ref{p1});}
\item[$(3)$]{The function
\begin{equation}\label{beta}
\beta(x)=\dfrac{V_d-V_c}{d-c}x + \dfrac{dV_c-cV_d}{d-c}
\end{equation}
is a $\mathcal{C}^2$-upper solution for the problem (\ref{p1}).}
\end{enumerate}
\end{proposition}

\noindent {\bf Proof.} \\
\begin{enumerate}

\item[$(1)$]{If $\dfrac{V_d}{d} \le \dfrac{V_c}{c}$, the thesis follows since 
$$
\alpha_1(c)=\dfrac{V_d}{d} c \le V_c, \quad \alpha_1(d)=V_d,$$ 
and, in $(c,d),$
$$ \alpha_1''(x) + H(x,\alpha_1(x),\alpha^\prime_1(x))=0.$$
Conversely, if $\alpha_1$ is a lower solution for problem (\ref{p1}), then $\alpha_1(c)=\dfrac{V_d}{d} c\leq V_c$, which implies $\dfrac{V_d}{d}\leq \dfrac{V_c}{c}.$}

\item[$(2)$]{Observe that $\alpha_k(c)=V_c$, $\alpha_k(d)=V_d$. On the other hand, as 
$$
H(x,\alpha_k(x),\alpha^\prime_2(x))\ge \dfrac{p(x)x^2 - \sqrt{p^2(x) x^4}- \sqrt{ 4x^3 q(x)|xz-y|}}{2x^3} 
$$ 
$$\ge -\sqrt{\dfrac{q(x)}{x^3}}\sqrt{|x \alpha^\prime_2(x)-\alpha_k(x)|},$$
$$
 \ge - \sqrt{\dfrac{Q}{c^3}}\sqrt{\max_{x \in [c,d]} \left|\dfrac{k}{2}(x^2-cd)+ \dfrac{V_d c-V_c d}{d-c}\right|},$$
then 
$$
\alpha^{\prime\prime}_2(x) + H(x,\alpha(x),\alpha^\prime_2(x)) \ge k- \sqrt{\dfrac{Q}{c^3}}\sqrt{\max_{x \in [c,d]} \left|\dfrac{k}{2}(x^2-cd)- V_c + \dfrac{V_d-V_c}{d-c}c\right|}.$$ 

Since $k >0$ satisfies the hiphothesis (\ref{k}), we derive
$$
\alpha_k''(x)+H(x,\alpha_k(x),\alpha^\prime_2(x))\ge 0.
$$
Then, the assertion holds.}

\item[$(3)$]{The thesis follows easily from the fact that 
$$\beta''(x)+H(x,\beta(x),\beta'(x))=H(x,\beta(x),\beta'(x)) \le 0$$ and
$$
\beta(c)=V_c, \quad \beta(d)=V_d.
$$}
\end{enumerate} \qed

\begin{remark} 1. Notice that there is no ambiguity in considering $k$ big enough such that the condition (\ref{k}) holds. In fact, it is easy to see that the maximum in (\ref{k}) depends on $k$ and is attained in the following way:
$$\max_{x \in [c,d]} \left|\dfrac{k}{2}(x^2-cd)+  \dfrac{V_d c-V_c d}{d-c}\right|=\max\left\{\left|\dfrac{k}{2}(d^2-cd)+  \dfrac{V_d c-V_c d}{d-c}\right|,\left|\dfrac{k}{2}(c^2-cd)+  \dfrac{V_d c-V_c d}{d-c}\right|\right\},$$
\noindent So, it is clear that in the proof of the previous proposition, we could choose  $k$ satisfying (\ref{k}) since
$$
\lim_{k \to \infty} \dfrac{k}{\sqrt{k}} = +\infty.$$
\smallskip

2. Observe also that the lower solution $\alpha_k$ can be written as $\alpha_k(x)= \beta (x) + \theta (x)$ where 
$$
\theta (x)= \dfrac{k}{2} \left(x^2 - (d+c)x +  cd \right). 
$$
In fact, 
$$
\alpha_k(x)=\dfrac{k}{2} x^2 + \left(\dfrac{V_d-V_c}{d-c}-\dfrac{k}{2}(d+c)\right)x + \dfrac{k}{2} cd - \dfrac{V_d c-V_c d}{d-c}$$ $$ 
\hspace{1,8cm}  =\left(\dfrac{V_d-V_c}{d-c}\right)x + \dfrac{dV_c-cV_d }{d-c}+\dfrac{k}{2} \left(x^2 - (d+c)x +  cd \right).
$$
The function $\theta$ is  quadratic, vanishes at $x=c$ and $x=d$ and is negative in $]c,d[$.\end{remark}

\begin{notation} Given two functions $\phi\leq \psi$ in $[c,d]$, let us denote by $[\phi,\psi]$ the functional interval
$$
[\phi,\psi]=\{ V \in W^{2,1}([c,d]) \, : \, \phi(x) \le V(x) \le \psi(x), \ \mbox{ for all $x \in [c,d]$}\}.$$
\end{notation}
Next, we state an existence and localization result for the problem (\ref{p1}).

\begin{theorem}\label{main1} Let $\alpha_1$, $\alpha_k$ and $\beta$ be the functions defined in the previous proposition.
\smallskip 

(a) If $\dfrac{V_d}{d} \le \dfrac{V_c}{c}$ then the problem (\ref{p1}) has extremal $W^{2,1}-$solutions, that is, the least and the greatest one, in the functional interval
$
[\alpha_1,\beta]$.
\smallskip 

(b) If $k >0$ satisfies (\ref{k}), then the problem (\ref{p1}) has extremal $W^{2,1}-$solutions in the functional interval $[\alpha_k,\beta].$
\smallskip

\end{theorem}

\noindent {\bf Proof.} Consider in case $(a)$
$$
E_1=\{(x,y,z) \in [c,d] \times \R^2 \, : \, \alpha_1(x) \le y \le \beta(x)\},$$
and in case $(b)$

$$
E_2=\{(x,y,z) \in [c,d] \times \R^2 \, : \, \alpha_k(x) \le y \le \beta(x)\}.$$
It is clear that $H$ is continuous in $E_1$ and $E_2$. Moreover, we have that
$$
|H(x,y,z)| \le \dfrac{p(x)}{x} + \sqrt{\dfrac{q(x)}{x^3} \beta(x)} + \dfrac{\sqrt{q(x)}}{x} \sqrt{|z|},$$
for all $(x,y,z) \in E_1$ or $(x,y,z) \in E_2$. Then, putting
\begin{equation}\label{nag1}
\hat{A}= \max_{x \in [c,d]} \left (\dfrac{p(x)}{x} + \sqrt{\dfrac{q(x)}{x^3} \beta(x)}\right ), \quad \hat{B}= \max_{x \in [c,d]} \dfrac{\sqrt{q(x)}}{x}\end{equation}
we derive that

\begin{equation}\label{nag2}
|H(x,y,z)| \le \hat{A} + \hat{B} \sqrt{|z|}\end{equation}

for $(x,y,z) \in E_1$ or $(x,y,z) \in E_2$, respectively. This inequality guarantees that the function $H$ satisfies the (classical) Nagumo condition both in $E_1$ and $E_2$. Using the fact that $\alpha_1,\alpha_k$ are $\mathcal{C}^2$-lower solutions in cases $(a)$ and $(b)$, respectively, and $\beta$ is a $\mathcal{C}^2$-upper solution for problem (\ref{p1}), such that 
$$ \alpha_1\leq \beta,\quad \alpha_k\leq \beta$$
the conclusion holds by application of a well--known result contained in \cite{cch}.\qed

\begin{corollary}\label{cor} Let $\dfrac{V_d}{d} \le \dfrac{V_c}{c}$ and $k >0$ satisfy (\ref{k}) and consider the functions $\alpha_1$, $\alpha_k$ and $\beta$ defined, respectively, by (\ref{alpha1}), (\ref{alpha2}), (\ref{beta}). Define 
$$\alpha(x)=\max\{\alpha_1(x),\alpha_k(x)\}.$$
Then the problem (\ref{p1}) has the extremal $W^{2,1}-$solutions in the interval $$[\alpha,\beta]=\{ V \in W^{2,1}([c,d]) \, : \, \alpha(x) \le V(x) \le \beta(x), \ \mbox{ for all $x \in [c,d]$}\}.$$

\end{corollary}

\noindent {\bf Proof.}  We observe that $$
\alpha(x)=\max\{\alpha_1(x),\alpha_k(x)\}$$ is a lower solution for the problem (\ref{p1}) in the sense defined in the previous section (not necessarily $\mathcal{C}^2$) and, as before, $\beta$ is a $\mathcal{C}^2$-upper solution for the problem (\ref{p1}). Similarly to the proof  of the  Theorem \ref{main1}, Nagumo condition holds in
$$E=\{(x,y,z) \in [c,d] \times \R^2 \, : \, \alpha(x) \le y \le \beta(x)\}.$$ Then, the result follows from an existence theorem contained in \cite{cch}. \qed

\begin{proposition}\label{Prop} The following assertions hold:

(a) Every solution $V$ of the problem (\ref{p1})  is convex. 
\smallskip

(b) Every solution $V$ of the problem (\ref{p1}) such that $V\geq \alpha_1$  satisfies    for all $x \in [c,d]$
$$xV'(x)-V(x) \le 0.$$
\end{proposition}

\noindent {\bf Proof.} \\
$(a)$ Clearly, the convexity of solutions of (\ref{p1}) derives from the fact that
$$
V''(x)=-H(x,V(x),V'(x) \ge 0 \ \mbox{ for all $x \in ]c,d[$}$$
and from the continuity in $c$ and $d$.\smallskip

\noindent$(b)$ Let $V$ be a solution of the problem (\ref{p1}) such that $V\geq \alpha_1$.  We claim that $\frac{V_d}{d}\geq V'(d)$. In fact, $\alpha_1=\frac{V_d}{d}x\leq V(x)$ implies that $\frac{V_d}{d}\geq \frac{V(x)-V_d}{x-d}$ in $[c,d]$ and then, letting $x\rightarrow d$, we obtain $\frac{V_d}{d}\geq V'(d)$. This inequality together with the fact that  the function $x \in [c,d] \longmapsto xV'(x)-V(x)$ is non-decreasing implies that $xV'(x)-V(x) \le 0$ for all $x \in [c,d]$.\qed  


Next theorem establishes the relations between convex solutions of (\ref{p1}) and of (\ref{p2})

\begin{theorem}\label{conv} Consider the problems (\ref{p1}) and  (\ref{p2}). Then

(a) If $\frac{V_d}{d}\leq \frac{V_c}{c}$, every solution $V$  of (\ref{p1}) provided by  the Theorem \ref{main1} is a convex solution of (\ref{p2}).

(b) Every convex solution  $V$ of the problem (\ref{p2}) is a convex solution of the problem (\ref{p1}).
\end{theorem}
\noindent {\bf Proof.} \\
As for (a), let $\frac{V_d}{d}\leq \frac{V_c}{c}$. Then $\alpha_1$ is a lower solution of (\ref{p1}) and every solution $V$  of (\ref{p1}) provided by  the Theorem \ref{main1} satisfies $V(x) \geq \alpha_1(x)$. Hence, by Proposition \ref{Prop}, $V$ is convex and 
 $$|x{V}'(x)-{V}(x)|=-(xV'(x)-V(x)).$$ 
 So
$$
{V}''(x)=-H(x,{V}(x),{V}'(x))= \dfrac{-p(x)x^2+\sqrt{p(x)x^4-4q(x)x^3(x V'(x)-V(x))}}{2x^3},$$
which shows clearly that ${V}$ is a convex function that solves (\ref{p2}). 

\noindent As for (b), let $V$ be a convex solution of (\ref{p2}). Then 
$$
x^3 (V'')^2(x) + p(x) x^2 V''(x) + q(x) (xV'(x)-V(x))=0,$$
which shows that 
$$ q(x) (xV'(x)-V(x))\leq 0,$$
since $V''\geq 0$, $p(x)\geq 0$ and $0<c\leq x \leq d.$
Therefore
$$
x^3 (V'')^2(x) + p(x) x^2 V''(x)+\frac{1}{4}p^2(x) x^4=\frac{1}{4}p^2(x) x^4- q(x) (xV'(x)-V(x)),$$
that is,
$$
(x^3 V''(x) + \frac{1}{2}p(x) x^2)^2=\frac{1}{4}p^2(x) x^4- q(x) (xV'(x)-V(x)).
$$
Applying the square root to both members, we obtain
$$
x^3 V''(x) + \frac{1}{2}p(x) x^2=\sqrt{\frac{1}{4}p^2(x) x^4- q(x) (xV'(x)-V(x))}.
$$
Then
$$
V''(x) =\frac{- p(x) x^2+\sqrt{p^2(x) x^4+| q(x) (xV'(x)-V(x))}}{2x^3 }=-H(x,{V}(x),{V}'(x)).
$$
So, ${V}$ solves (\ref{p2}). \qed

From Theorem \ref{main1}, Theorem \ref{conv} and Corollary \ref{cor}, it is clear that the following existence and localization result holds:

\begin{theorem}\label{dir} Consider the problem (\ref{p}) with standard Dirichlet conditions, that is,
\begin{equation*}
\left\{
\begin{array}{ll}
x^3 (V'')^2(x) + p(x) x^2 V''(x) + q(x) (xV'(x)-V(x))=0, \\
\\
V(c)=V_c, \quad V(d)=V_d.
\end{array}
\right.
\end{equation*}

\begin{enumerate}

\item[$(1)$] {If $\dfrac{V_d}{d} \le \dfrac{V_c}{c}$ then this problem has the extremal convex $W^{2,1}-$solutions in the functional interval $[\alpha_1,\beta]$, where $\alpha_1$ and $\beta$ are provided, respectively, by (\ref{alpha1}) and (\ref{beta});}

\item[$(2)$] {If If $\dfrac{V_d}{d} \le \dfrac{V_c}{c}$ and $k>0$ satisfies (\ref{k}) then this problem has the extremal convex $W^{2,1}$-solutions in the functional interval $[\alpha,\beta]$, where $$\alpha(x)=\max\{\alpha_1(x),\alpha_k(x)\}$$ and $\alpha_1$, $\alpha_k$, $\beta$ are provided, respectively, by (\ref{alpha1}), (\ref{alpha2}), (\ref{beta}).} 
\end{enumerate}
\end{theorem}

\begin{remark} Under the hyphoteses of the above theorem, observe that if $\dfrac{V_d}{d}=\dfrac{V_c}{c}$ then $\alpha_1$ is a solution of the problem (\ref{p2}). On the other hand, in the periodic case, $V_c=V_d$, the constant function $V\equiv V_c$ is a solution of (\ref{p2}).
\end{remark}

\begin{example}\label{ex1}
Consider problem (\ref{p2}) in the interval $[c,d]=[2,6]$, with $$p(x)=1+x^3, \quad q(x)=[x],$$ where $[ \cdot ]$ denotes integer part,  and boundary conditions $V(2)=9$, $V(6)=1$. Notice that in this case it is $V_c > V_d$. Condition (\ref{k}) says that we have to get $k >0$ such that
$$
k \ge \sqrt{\dfrac{6}{8}} \sqrt{\max \{|-4k-13|,|12k-13|\}},$$

so simple computations show that $k=8$ satisfies (\ref{k}). Then, function
$$
\alpha(x)=4x^2 - 34x + 61$$
is a lower solution for this problem. On the other hand,
$$
\beta(x)=13-2x$$
is an upper solution. \\

Then, by application of Theorem \ref{dir}, problem (\ref{p2}) with $p(x)=1+x^3$, $q(x)=[x]$, and boundary conditions $V(2)=9$, $V(6)=1$, has the extremal solutions in the functional interval
$$
\left[\max\left\{\dfrac{1}{6}x,4x^2-34x+61\right\}, 13-2x\right].$$

\end{example}

\section{Problem with functional boundary conditions}

In this section we deal with problem (\ref{pp}) on its full expression and, as said in introduction, we will use a generalized monotone method. In the construction of this method
we will use two technical lemmas. First of them is the following generalization of Bolzano's theorem.

\begin{lemma}\label{bolzano}{\bf \cite[Lemma 2.3]{danipouso}}
Let $a,b \in \mathbb R$,
$a \leq b$, and
$h:\mathbb R \longrightarrow \mathbb R$ a function satisfying $h(a) \leq 0
\leq h(b)$
and
\begin{equation}
\label{bol}
\liminf_{z \to x^-}h(z) \geq h(x) \geq \limsup_{z \to x^+}h(z)
\
\mbox{ for all $x \in [a,b]$.}
\end{equation}
Then there exist $c_1,c_2 \in [a,b]$ such that $h(c_1)=0=h(c_2)$
and if $h(c)=0$ for some $c\in [a,b]$ then $c_1 \leq
c \leq c_2$,
that is, $c_1$ and $c_2$ are, respectively, the least zero and the greatest one of $h$ in $[a,b]$.
\end{lemma}

The second auxiliar result we need deals with the existence of extremal fixed points for nondecreasing operators defined in the space of absolutely continuous functions.

\begin{lemma}\label{heikkila}{\bf \cite[Proposition 1.4.4]{hela}} Let $I \subset \R$ a nonempty closed interval and $[\alpha,\beta]$ a nonempty functional interval in $AC(I)$. Assume that $G:[\alpha,\beta] \longrightarrow [\alpha,\beta]$ is a nondecreasing mapping and that there exists $\psi \in L^1(I,[0,+\infty))$ such that
$$
|(G\gamma)'(x)| \le \psi(x) \ \mbox{for all $\gamma \in [\alpha,\beta]$ and almost all $x \in I$}.$$
Then, $G$ has in $[\alpha,\beta]$ the greatest, $V^*$, and the least, $V_*$, fixed points. Moreover, they satisfy
\begin{equation}\label{maxmin}
V_*=\min\{V \, : \, GV \le V\}, \quad V^*=\max\{V \, : \, V \le GV\}.\end{equation}
\end{lemma}

Now we establish a new result on the existence of extremal convex solutions for problem (\ref{pp})

\begin{theorem}\label{main3} Assume that there exist $\alpha,\beta \in AC([c,d])$ which are, respectively, a lower and an upper solution for problem (\ref{pp}) satisfying $\alpha(x) \le \beta(x)$ for all $x \in [c,d]$. Put
$$
[\alpha,\beta]=\{\gamma \in AC([c,d]) \, : \, \alpha(x) \le \gamma(x) \le \beta(x) \ \mbox{ for all $x \in [c,d]$}\},$$
$$
E=\left\{y \in \R \, : \, \min_{x \in [c,d]} \alpha(x) \le y \le \max_{x \in [c,d]} \beta(x)\right\},$$
and assume, moreover, that the following conditions hold:

\begin{enumerate}

\item[$(H_1)$] For all $\gamma \in [\alpha,\beta]$ and all $y \in E$ we have
$$
\liminf_{z \to y^-} B_i(z,\gamma) \ge B_i(y,\gamma) \ge \limsup_{z \to y^+} B_i(z,\gamma) \ \ (i=1,2);$$

\item[$(H_2)$] For all $y \in E$ the functions $B_i(y,\cdot)$ are nonincreasing in $[\alpha,\beta]$ $(i=1,2)$, that is, if $\gamma_1,\gamma_2 \in [\alpha,\beta]$ are such that $\gamma_1(x) \le \gamma_2(x)$ for all $x \in [c,d]$, then $B_i(y,\gamma_1) \ge B_i(y,\gamma_2)$.

\end{enumerate}

In these conditions, problem (\ref{pp}) has the extremal convex solutions in $[\alpha,\beta]$.

\end{theorem}

\noindent {\bf Proof.} \\
We define a mapping $G:[\alpha,\beta] \longrightarrow [\alpha,\beta]$ as follows: for all $\gamma \in [\alpha,\beta]$, $G\gamma$ is the greatest convex solution in $[\alpha,\beta]$ for the Dirichlet problem
\begin{equation}\label{aux}
\left\{
\begin{array}{ll}
V''(x)+H(x,V(x),V'(x))=0 \ \mbox{ for all $x \in [c,d]$,} \\
\\
V(c)=\gamma_c, \quad V(d)=\gamma_d,
\end{array}
\right.
\end{equation}
where $\gamma_c,\gamma_d$ are the greatest solutions in, respectively, $[\alpha(c),\beta(c)]$ and $[\alpha(d),\beta(d)]$, for the following respective algebraic equations:
\begin{eqnarray}
\label{c}
B_1(y,\gamma)=0, \\
\label{d}
B_2(y,\gamma)=0.
\end{eqnarray}

{\it Step $1$: The mapping $G$ is well--defined.} First, by virtue of being $\alpha$ and $\beta$ a lower and an upper solution for problem (\ref{pp}) and by condition $(H_2)$, we have for all $\gamma \in [\alpha,\beta]$:
$$
\begin{array}{cc}
B_1(\alpha(c),\gamma) \le B_1(\alpha(c),\alpha) \le 0 \le B_1(\beta(c),\beta) \le B_1(\beta(c),\gamma), \\
\\
B_2(\alpha(d),\gamma) \le B_2(\alpha(d),\alpha) \le 0 \le B_2(\beta(d),\beta) \le B_2(\beta(d),\gamma)
\end{array}
$$
and so condition $(H_1)$ implies that the numbers $\gamma_c$ and $\gamma_d$ are well--defined, by application of Lemma \ref{bolzano}. \\

On the other hand, the fact of being $\gamma_c$ and $\gamma_d$ the greatest solutions of equations (\ref{c})--(\ref{d}) in, respectively, $[\alpha(c),\beta(c)]$ and $[\alpha(d),\beta(d)]$, implies that $\alpha(c) \le \gamma_c$, $\alpha(d) \le \gamma_d$, $\beta(c) \ge \gamma_c$ and $\beta(d) \ge \gamma_d$. So, $\alpha$ and $\beta$ are, respectively, a lower and an upper solution for problem (\ref{aux}). This guarantees that (\ref{aux}) has the greatest convex solution in $[\alpha,\beta]$. (Notice that (\ref{nag1})--(\ref{nag2}) provides a Nagumo--type bound for $H$ between our $\alpha$ and $\beta$.) \\

{\it Step $2$: $G$ is a nondecreasing mapping.} Let $\gamma_1,\gamma_2 \in [\alpha,\beta]$ such that $\gamma_1(x) \le \gamma_2(x)$ for all $x \in [c,d]$ and we will show that $G\gamma_1 \le G\gamma_2$. First, notice that
$$
B_1(G\gamma_1(c),\gamma_2) \le B_1(G\gamma_1(c),\gamma_1)=0$$
and
$$
B_1(\beta(c),\gamma_2) \ge B_1(\beta(c),\beta) \ge 0,$$
so reasoning as above we obtain that $G\gamma_1(c) \le G\gamma_2(c)$. In a similar way we prove that $G\gamma_1(d) \le G\gamma_2(d)$. Now, assume that $G\gamma_1 \nleqslant G\gamma_2$ and consider the function
$$
\hat{\alpha}(x)=\max\{G\gamma_1(x),G\gamma_2(x)\}.$$
Thus defined, $\hat{\alpha}$ is a lower solution for problem (\ref{aux}) with conditions $V(c)=\gamma_{2c}$, $V(d)=\gamma_{2d}$. Indeed, if $G\gamma_1 < G\gamma_2$ in an interval $(\hat{x}_1,\hat{x}_2)$ then
$$
\hat{\alpha}''(x)=G\gamma_2''(x)=-H(x,G\gamma_2(x),(G\gamma_2)'(x))=-H(x,\hat{\alpha}(x),\hat{\alpha}'(x)) \ \mbox{ for all $x \in (\hat{x}_1,\hat{x}_2)$}.$$ In the same way, if $G\gamma_1 > G\gamma_2$ in an interval $(\tilde{x}_1,\tilde{x}_2)$ then
$$
\hat{\alpha}''(x)=-H(x,\hat{\alpha}(x),\hat{\alpha}'(x)) \ \mbox{ for all $x \in (\tilde{x}_1,\tilde{x}_2)$}.$$
On the other hand, if $G\gamma_1(x_0)=G\gamma_2(x_0)$ then convexity of $G\gamma_1$ and $G\gamma_2$ implies that $D_-\hat{\alpha}(x_0) < D^+\hat{\alpha}(x_0)$. \\
So, problem (\ref{aux}) with conditions $V(c)=\gamma_{2c}$, $V(d)=\gamma_{2d}$, has a solution in $[\hat{\alpha},\beta]$, but this contradicts the fact that $G\gamma_2$ is the greatest solution for this problem in $[\alpha,\beta]$. Then, we conclude that $G\gamma_1 \le G\gamma_2$ and so $G$ is a nondecreasing mapping. \\

{\it Step $3$: $G$ has the extremal fixed points.} Let $\gamma \in [\alpha,\beta]$. For all $x \in [c,d]$ we have that
$$
(G\gamma)'(x)=(G\gamma)'(c)- \int_c^x H(s,G\gamma(s),G\gamma'(s)) \, ds,$$
and so
$$
|G\gamma'(x)| \le \max\left\{ \dfrac{|\beta(d)-\alpha(c)|}{d-c}, \dfrac{|\beta(c)-\alpha(d)|}{d-c}\right\} + \int_c^x |\hat{A} + \hat{B} \sqrt{s}| \, ds,$$
where $\hat{A}$, $\hat{B}$ are as in (\ref{nag1}). \\
Then, by application of Lemma \ref{heikkila} we obtain that $G$ has the extremal fixed points in $[\alpha,\beta]$, say $V^*$, $V_*$, which moreover satisfy (\ref{maxmin}). \\

{\it Step $4$: $V^*$ is the greatest convex solution of problem (\ref{pp}) in $[\alpha,\beta]$.} First, it is clear, as $GV^*=V^*$, that $V^*$ is a solution of problem (\ref{pp}). Now, if $V$ is another solution of (\ref{pp}) then we have that $V \le GV$ and so (\ref{maxmin}) implies that $V \le V^*$. So, $V^*$ is the greatest convex solution of problem (\ref{pp}) in $[\alpha,\beta]$. \\

To obtain the least convex solution of (\ref{pp}) in $[\alpha,\beta]$ we only have to redefine the mapping $G$ in the obvious way. \qed

\begin{remark} Notice that condition $(H_1)$ is satisfied, for example, if $B_1(\cdot,\gamma)$ is continuous or if it has only downwards discontinuities.
\end{remark}

\begin{remark} The same argument used in the proof of Theorem \ref{dir} provides that extremal convex solutions of problem (\ref{pp}) are extremal convex solutions of problem (\ref{p}) too.
\end{remark}

\begin{example}\label{exfun} Consider problem (\ref{p}) in an interval $[c,d]$, $c>0$, with the following boundary conditions:

\begin{enumerate}

\item[$(B_1)$] ``The initial value of the solution is one half of its mean value on the whole interval $[c,d]$'';

\item[$(B_2)$] ``The final value of the solution has integer part $4$'';

\end{enumerate}

Previous conditions can be written in the following form:

\begin{enumerate}

\item[$(B_1)$] $B_1(V(c),V)=V(c)-\dfrac{1}{2} \dfrac{1}{d-c} \displaystyle{\int_c^d V(s) \, ds}=0$;

\item[$(B_2)$] $B_2(V(d),V)=-[V(d)]+4=0$.

\end{enumerate}

We will show that for $V(d) \in [4,5)$ and $d \ge 3c$, $\alpha(x)=\dfrac{V(d)}{d}x$ and $\beta \equiv V(d)$ are, respectively, a lower and an upper solution for this problem. Indeed, we have $-[\alpha(d)]=-4$ and
$$
\alpha(c)-\dfrac{1}{2} \dfrac{1}{d-c} \int_c^d \alpha(s) \, ds= \dfrac{V(d)}{d} \left(c-\dfrac{d+c}{4}\right) \le 0$$
if $d \ge 3c$. On the other hand, the constant function $\beta \equiv V(d)$ is such that
$$
B_1(\beta(c),\beta)=\dfrac{1}{2} V(d) \ge 0$$
and
$$
B_2(\beta(d),\beta)=-[V(d)]+4=0.$$
Finally, notice that for each $V(d) \in [4,5)$, functions $B_1$ and $B_2$ satisfy conditions $(H_1)-(H_2)$ between $\alpha$ and $\beta$. We can conclude, by application of Theorem \ref{main3}, that if $d \ge 3c$ then problem (\ref{pp}) with boundary conditions $(B_1)-(B_2)$ has the extremal solutions between $\alpha$ and $\beta$ for each $V(d) \in [4,5).$

\end{example}

\begin{remark} In papers \cite{amsab} and \cite{gromor} the authors obtained uniqueness of solutions for problem (\ref{dir}). Notice that uniqueness of solutions cannot be guaranteed when we include variable coefficients $p$ and $q$ and functional boundary conditions. For example, consider problem (\ref{p}) with boundary conditions
$$
[V(1)]=2, \quad [V(2)]=4.$$
In this case, each function $V(x)=ax$ with $a \in [2,2.5)$ is a solution of the problem.
\end{remark}

\textbf{Acknowledgment}: The first author was partially supported by Xunta de Galicia, Conseller\'ia de Cultura, Educaci\'on e Ordenaci\'on Universitaria, through the project EM2014/032 ``Ecuaci\'ons diferenciais non lineares"; and by by Ministerio de Econom\'ia y Competitividad of Spain under Grant MTM2010-15314, cofinanced by the European Community fund FEDER. \\
The second author was partially funded by Funda\c c\~ao para a Ci\^encia e Tecnologia
through the project UID/Multi/00491/2013 and the Transnational Cooperation FCT Portugal-Slovakia ``Analysis of Nonlinear Partial Differential Equations in Mathematical Finance (2013-2014)'' and by the EU Grant Program FP7-PEOPLE-2012-ITN STRIKE - ``Novel
Methods in Computational Finance'', No. 304617 (D.S.).

\end{document}